\documentclass[twoside,12pt]{article}
\usepackage{geometry}
\usepackage{amsmath,amssymb,amsfonts,amsthm}
\usepackage{pxfonts}
\usepackage{txfonts}
\usepackage{todonotes}
\usepackage{tablists}
\usepackage{graphicx}
\usepackage{mathrsfs}
\usepackage{color}
\usepackage{float}
\usepackage{extarrows}
\usepackage{exscale}
\usepackage{titlesec}
\titleformat{\section}[block]{\large\center\sc}{\arabic{section}}{0.5em}{}[] % ����section ��ʽ
\usepackage{relsize}
\usepackage{color}
\definecolor{teal}{RGB}{0,128,172}
\definecolor{pur}{RGB}{224,104,255}
\usepackage{hyperref}
\hypersetup{hypertex=true,
            colorlinks=true,
            linkcolor=blue,
            anchorcolor=teal,
            citecolor=pur}
\usepackage[titles]{tocloft}
\usepackage{xcolor}
\definecolor{teal}{RGB}{67,205,128}
\definecolor{SBlue}{RGB}{136,120,255}

\usepackage{fancyhdr}
\theoremstyle{plain}
\newtheorem{theorem}{Theorem}[section]
\newtheorem{lemma}[theorem]{Lemma}

\newtheorem{proposition}[theorem]{Proposition}

\newtheorem{remark}[theorem]{Remark}

\let\oldsection\section
\renewcommand\section{\setcounter{equation}{0}\oldsection}

\def\be{\begin{equation}}
\def\ee{\end{equation}}
\def\bes{\begin{equation*}}
\def\ees{\end{equation*}}
\def\bs{\begin{split}}
\def\es{\end{split}}
\def\bali{\begin{aligned}}
\def\eali{\end{aligned}}
\newcommand{\pf}{\noindent {\bf Proof. \hspace{2mm}}}
\def\bR{{\mathbb R}}
\def\un{\underbrace}
\def\al{\alpha}

\def\th{\theta}

\def\Dl{\Delta}
\def\lt{\left}
\def\rt{\right}

\def\i{\infty}
\def\p{\partial}
\def\f{\frac}
\def\na{\nabla}

\def\O{\Omega}

\def\q{\quad}

\def\bl{\boldsymbol}

\def\mR{\mathbb{R}}

\allowdisplaybreaks
\def\mD{\mathcal{D}}

\begin{document}

\title{\Large\bf A refined uniqueness result of Leray's problem in an infinite-long pipe with the Navier-slip boundary}

%On Leray's problem in an infinite-long pipe with the Navier-slip boundary condition
%\author{Zijin Li $^{a,}$\footnote{E-mail:zijinli@smail.nju.edu.cn}\vspace{0.5cm}\\
%\footnotesize $^a$School of Mathematics and Statistics, Nanjing University of Information Science $\&$ Technology,\\
%\footnotesize Nanjing 210044, China.\\
%\vspace{0.5cm}
%}
\author{\normalsize\sc {Zijin Li, Ning Liu and Taoran Zhou}}

\date{}

\maketitle
% \vskip 0.2in

\begin{abstract}
We consider the generalized Leray's problem with the Navier-slip boundary condition in an infinite pipe $\mD=\Sigma\times\mR$. We show that if the flux $\Phi$ of the solution is no larger than a critical value that is \emph{independent with the friction ratio} of the Navier-slip boundary condition, the solution to the problem must be the parallel Poiseuille flow with the given flux. Compared with known related 3D results, this seems to be the first conclusion with the size of critical flux being uniform with the friction ratio $\al\in]0,\infty]$, and it is surprising since the prescribed uniqueness breaks down immediately when $\al=0$, even if $\Phi=0$.

 Our proof relies primarily on a refined gradient estimate of the Poiseuille flow with the Navier-slip boundary condition. Additionally, we prove the critical flux $\Phi_0\geq\f{\pi}{16}$ provided that $\Sigma$ is a unit disk.

\medskip


\end{abstract}

\tableofcontents

\section{Introduction}
\q\ The 3D stationary Navier-Stokes (NS) equations read:
\be\label{NS}
\lt\{
\begin{aligned}
&\bl{u}\cdot\na \bl{u}+\na \mathcal{\pi}-\Dl \bl{u}=0,\\
&\na\cdot \bl{u}=0,
\end{aligned}
\rt.\q \text{in}\q \mathcal{D}\subset \bR^3.
\ee
Here  $\bl{u}(x)\in\mathbb{R}^3$, $\mathcal{\pi}(x)\in\mathbb{R}$ represent the velocity and the scalar pressure respectively. The domain $\mD$ is an infinite long pipe: $\mD=\Sigma\times\mR$, where $\Sigma$ is a bounded smooth 2D domain in the $x_h=(x_1,x_2)$ direction.
\begin{figure}[H]\label{FIG1}
\centering
\includegraphics[scale=0.6]{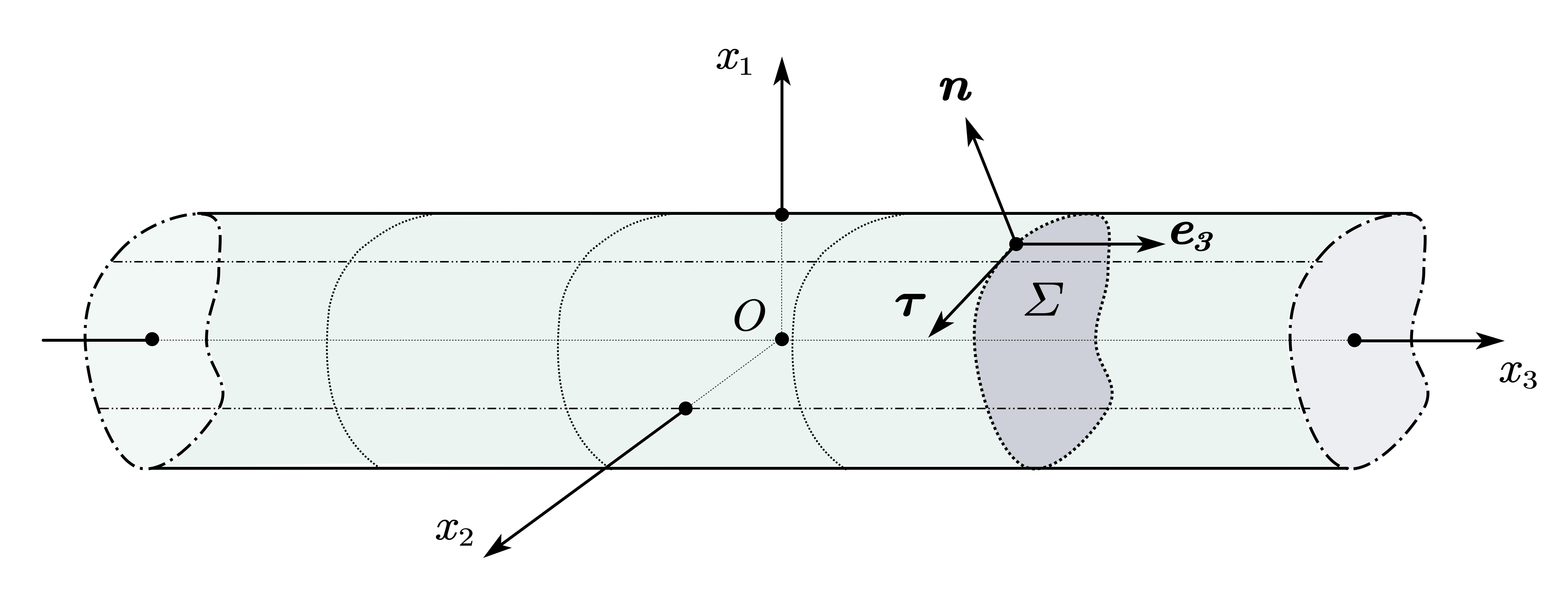}
\caption{The infinite long pipe $\mD=\Sigma\times\mR$.}
\end{figure}
We impose the Navier-slip boundary condition to the velocity field $\bl{u}$ on $\p\mD$:
\be\label{NBC}
\left\{
\begin{aligned}
&2(\mathbb{S}\bl{u}\cdot\boldsymbol{n})_{\mathrm{tan}}+\alpha \bl{u}_{\mathrm{tan}}=0,\\
&\bl{u}\cdot\boldsymbol{n}=0,\\
\end{aligned}
\right.\q\text{on}\q\p\mD.
\ee
Here $\mathbb{S}\bl{u}=\frac{1}{2}\left(\nabla \bl{u}+\nabla^T \bl{u}\right)$ is the stress tensor, where $\nabla^T \bl{u}$ stands for the transpose of the Jacobian matrix $\nabla \bl{u}$, and $\boldsymbol{n}$ is the unit outer normal vector of $\p \mD$. For a vector field $\bl{v}$, its tangential part $\bl{v}_{\mathrm{tan}}$ is defined by
\[
\bl{v}_{\mathrm{tan}}:=\bl{v}-(\bl{v}\cdot\boldsymbol{n})\boldsymbol{n}.
\]
The constant $\al>0$ stands for the friction ratio, which may depend on the property of the boundary,  the viscosity of the fluid, et cetera. We mention that when $\al\to0_+$, the boundary condition \eqref{NBC} turns to be the total Navier-slip boundary condition, while when $\al\to+\infty$, the boundary condition \eqref{NBC} degenerates into the no-slip boundary condition $\bl{u}\big|_{\p\mD}\equiv 0$. In this paper, we consider the intermediate case $0<\al<+\i$.

Let $(\bl{u},p)$ be a smooth solution of the boundary value problem \eqref{NS}--\eqref{NBC}, we define
\be\label{Fluxc}
\Phi=\int_{\Sigma}u_3(x_h,x_3)dx_h,
\ee
the \emph{flux} of the flow. In this paper, we suppose $\Phi\geq 0$ without loss of generality. We mention that $\Phi$ is independent of $x_3$, owing to the divergence-free property and the impermeable boundary condition of $\bl{u}$:
\[
\begin{split}
&\f{d}{dx_3}\int_{\Sigma}u_3(x_h,x_3)dx_h=\int_{\Sigma}\p_{x_3}u_3(x_h,x_3)dx_h\\
=&-\int_{\Sigma}\left(\p_{x_1}u_1(x_h,x_3)+\p_{x_2}u_2(x_h,x_3)\right)dx_h=-\int_{\p\Sigma}\bl{u}\cdot\bl{n}ds=0\,.
\end{split}
\]

Given a flux $\Phi\geq0$, we say $\bl{u}$ is a weak solution of \eqref{NS}--\eqref{NBC}--\eqref{Fluxc}, \emph{the generalized Leray's problem with the Navier-slip boundary condition}, if and only if the following (i) and (ii) hold:
\begin{itemize}
\item[(i)]{$\bl{u}\in H^1_\sigma(\overline{\mD})$} satisfies \eqref{NS}--\eqref{NBC} weakly. That is, for any $\bl{\varphi}\in C_c^{\infty}\left(\overline{\mathcal{D}}\right)$ with $\na\cdot\bl{\varphi}=0$ and $\bl{\varphi}\cdot\bl{n}\Big|_{\partial\mathcal{D}}=0$, one has
\[
2 \int_{\mathcal{D}} \mathbb{S} \bl{u}: \mathbb{S}\bl{\varphi} d x+\alpha \int_{\partial \mathcal{D}} \boldsymbol{u}_{\tan } \cdot \bl{\varphi}_{\tan } d S-\int_{\mathcal{D}} \boldsymbol{u} \cdot \nabla \bl{\varphi} \cdot \boldsymbol{u} d x=0\,;
\]

\item[(ii)]$\bl{u}$ satisfies \eqref{Fluxc} in the sense of trace.
\end{itemize}

We denote $\bl{\tau}$ and $\bl{n}$ the unit tangent vector and the unit outer normal vector of the cross-section $\Sigma$, respectively. See Figure \ref{FIG1}. Noticing that the coordinate vector $\bl{e_3}$ is always a tangent vector on $\p\mD$, we can rewrite
\be\label{NaCo}
\bl{u}=u_\tau\bl{\tau}+u_3\bl{e_3}+u_n\bl{n}
\ee
near the boundary. Under this curvilinear coordinate system, the Navier-slip boundary condition \eqref{NBC} can be rewritten as the following Robin-Dirichlet-mixed form:
\be\label{NBCM1}
\left\{\begin{array}{l}
\partial_{\boldsymbol{n}} u_{\tau}=\left(\kappa(x)-\alpha\right) u_{\tau}, \\
\partial_{\boldsymbol{n}} u_3=-\alpha u_3, \\
u_n=0,
\end{array} \q\text { on } \partial \mathcal{D}.\right.
\ee
Here $\kappa$ is the curvature of $\p\Sigma$, which is uniformly bounded since $\Sigma$ is bounded and smooth. For the detailed derivation of \eqref{NBCM1}, we refer readers to \cite[Proposition 2.1 and Corollary 2.2]{Watanabe} (also \cite[Section 2.1]{LPY2024SCM}).

With the help of the representation \eqref{NBCM1}, it is clear that the following Poiseuille flow is an exact solution to the problem \eqref{NS}--\eqref{NBC}--\eqref{Fluxc}:
\be\label{HP1}
\left\{
\begin{aligned}
&\bl{P}_\Phi=P_\Phi\bl{e_3}, \\
&-\Dl_hP_\Phi(x_h)=Constant,\q\q\q\text{in }\Sigma,\\
&\f{\p P_\Phi}{\p\bar{\bl{n}}}=-\al P_\Phi,\hskip 1.9cm\q\q\q\text{on }\p\Sigma,\\
&\int_{\Sigma}P_\Phi(x_h)dx_h=\Phi. \\
\end{aligned}
\right.
\ee
Here $\bar{\bl{n}}$ is the unit outer normal vector of $\p\Sigma$, while $P_\Phi:\Sigma\to\mathbb{R}$ is a scaler function. Its related pressure field reads:
\[
\pi_P=-Constant\cdot x_3,\q\q\q\text{in }\Sigma\,.
\]
Now a problem arises:

\begin{center}
\fbox{\emph{Is \eqref{HP1} the only solution to the problem \eqref{NS}--\eqref{NBC}--\eqref{Fluxc}?}}
\end{center}
\subsection{Main result}
Our main result of the current paper states:
\begin{theorem}\label{th1}
{\sl Let $(\bl{u},p)$ be a weak solution of the boundary value problem \eqref{NS}--\eqref{NBC}, and
\[
\Phi=\int_{\Sigma}u_3(x_h,x_3)dx_h\,.
\]
Then there exists a critical flux $\Phi_0>0$, which is independent with the friction ratio $\al>0$, that if $\Phi\leq\Phi_0$, and
\be\label{GROWC}
\|\na\bl{u}\|_{L^2(\Sigma\times\,]-\zeta,\zeta[\,)}=o(\zeta^{3/2}),\q\text{for}\q\zeta>>1\,,
\ee
then $\bl{u}$ must be the Poiseuille flow $\bl{P}_\Phi$ given in \eqref{HP1}.}
\end{theorem}

\qed

\begin{remark}
{\sl Here the growth condition \eqref{GROWC} is relatively mild since we mainly focus on bounded solutions. However, the uniform bound of the critical flux $\Phi_0$ as $\al\to 0_+$ seems surprising. Indeed, other than the parallel flow \eqref{HP1}, some parasitic solutions may exist in the total slip case. For example, if the cross section $\Sigma$ is a unit disk, the helical solution:
\[
\bl{u_h}=Cr\bl{e_\th}+\f{\Phi}{\pi}\bl{e_3}\,,
\]
solves \eqref{NS}--\eqref{NBC}--\eqref{Fluxc} when $\al=0$ for any $C\in\mathbb{R}$. Here
\[
\bl{e_\th}=\Big(\f{-x_2}{\sqrt{x_1^2+x_2^2}},\f{x_1}{\sqrt{x_1^2+x_2^2}},0\Big)\,.
\]
This breaks down the desired uniqueness even when $\Phi=0$.}
\end{remark}

\qed

\begin{remark}
{\sl In Section \ref{Sec4.4}, we show the critical flux $\Phi_0\geq\f{\pi}{16}$ when $\Sigma$ is a unit disk.}
\end{remark}

\qed

\begin{remark}
{\sl Compared with the result in \cite{LPY2024SCM} which showed the uniqueness provided $\Phi$ is no bigger than a critical flux { $\Phi_0=C(\alpha, \mD)$}, Theorem \ref{th1} excludes the dependence of the friction ratio $\al$ for $\Phi_0$. Although the authors did not provide a specific expression  of $\Phi_0$ there, readers can verity $\Phi_0\simeq_\mD\f{\al}{1+\al}$ by applying the method  therein. Theorem \ref{th1} improves the result in \cite{LPY2024SCM} especially when $\al$ is small (close to the total slip Navier boundary condition).}
\end{remark}

\qed

\subsection{Related works}
We review some works related to the solvability of Leray's problem. Originally, the study of Leray's problem focuses on the existence, regularity, uniqueness, and asymptotic behavior of Systems \eqref{NS}--\eqref{Fluxc} subject to no-slip boundary condition
\[
\bl{u}=0,\q\text{on } \p\mD\,,
\]
as described by Ladyzhenskaya of her works \cite[p. 77]{MR0119676} and \cite[p. 551]{MR0106641}. Later in\cite{LaSo}, Ladyzhenskaya and Solonnikov provided a comprehensive analysis of the existence, uniqueness, and asymptotic behavior of weak solutions to the Leray's problem with the help of small-flux conditions. Amick \cite{MR0510120,MR0512162} made a significant initial contribution to solving Leray's problem by reducing the issue to a variational problem, although the uniqueness of the solutions remains an open question. Further details on the decay, far-field asymptotic behavior, and well-posedness of solutions to Leray's problem under the no-slip boundary condition are available in \cite{MR1000723,MR0502782,MR1992105}. For a systematic exploration and analysis of this problem, readers may consult Chapter XIII of \cite{Galdi2011}. In recent research, Wang and Xie investigated the existence, uniqueness, and uniform structural stability of Poiseuille flows in a pipe, addressing the 3D axially symmetric inhomogeneous Navier-Stokes equations. See \cite{MR4453236,MR4381142} and some of their subsequent works.

The Naiver-slip boundary condition \eqref{NBC}, which was initialed by Navier \cite{Navier}, allows fluid to slip along the boundary with a scale being proportional to its stress tensor.  Due to the complexity of the Navier-slip boundary condition, new mathematical strategies are needed to study problems related to it. More and more studies on problems with Navier-slip boundary conditions emerged these years, see \cite{AACG} and references therein. In \cite{Mucha2003, Mucha:2003STUDMATH, Konie:2006COLLMATH}, authors studied the solvability of the steady Navier-Stokes equations with the total Navier-slip condition ($\al=0$). Wang and Xie \cite{WangXie2023JDE} proved the uniqueness and uniform structural stability of the Poiseuille flow in an infinite straight long pipe with the Navier-slip boundary condition. Li, Pan and Yang  \cite{LPY2023Ar} gave the characterization of bounded smooth solutions to the axisymmetric Navier-Stokes equations with the total Navier-slip boundary condition and the Navier-Hodge-Lions boundary condition in the infinitely long cylinder.

In the recent paper \cite{LPY2024SCM}, authors proved the generalized Leray's problem of stationary 3D Navier-Stokes equations (i.e. existence, uniqueness, together with the higher-order regularity and asymptotic behavior of the stationary 3D Navier-Stokes flow in a distorted infinite pipe with Navier-slip boundary condition) with the flux of flow being small enough, depending on the shape of the pipe and the friction ratio $\al$ of the Naiver-slip boundary condition. Meanwhile, the related 2D problem, which is studied in \cite{LPY2023JDE}, presented an exact formula that measures the smallness of $\Phi$ with respect to $\al$ when considering the existence, regularity together with the asymptotic behavior problem:
\be\label{E1}
\Phi\leq\Phi_0:=\f{C(1+\al)}{\al}\,,
\ee
where $C>0$ is a constant independent with $\al$. Estimate \eqref{E1} shows when $\al$ is small enough, the smallness restriction of the flux tends weaker, and for the total Navier-slip condition ($\al=0$), the solution can be arbitrarily large.

\subsection{Difficulties, the outline of proof, and the structure of paper}
The cross section $\Sigma$ degenerates to an interval in the 2D case, and we assume $\Sigma=]0,1[$ without loss of generality. Direct calculation shows
\be\label{2dpoi}
\bl{P}_{\Phi}=\lt(0\,\,,\,\, \f{6\al\Phi}{6+\al}\lt(-x^2_1+x_1\rt)+\f{6\Phi}{6+\al}\rt),\q(x_1,x_2)\in\,]0,1[\,\times\mR\,,
\ee
solves \eqref{HP1} (with $\bl{e_3}$ replaced by $\bl{e_2}$). It is clear that
\be\label{E2}
|\nabla \bl{P}_\Phi|\leq\f{C\Phi\al}{1+\al}\,.
\ee
For the 3D case, the gradient bound \eqref{E2} also holds when the cross-section $\Sigma$ is a unit disk. Indeed,
\be\label{Hagen}
\bl{P}_{\text{disk},\Phi}(x):=\f{2(\al+2)\Phi}{(\al+4)\pi}\lt(1-\f{\al}{\al+2}|x_h|^2\rt)\boldsymbol{e_3},\q\text{with its pressure } \pi_\Phi(x)=-\f{8\al\Phi}{(\al+4)\pi} x_3,
\ee
which could be considered as a generalization of the Hagen-Poiseuille flow under the no-slip boundary condition ($\al=+\i$), solves the problem \eqref{HP1} provided $\Sigma=\{x_h\in\mR^2:\,|x_h|<1\}$. In this case, clearly
\be\label{E1777}
|\na\bl{P}_{\text{disk},\Phi}|\lesssim\f{\al\Phi}{1+\al}
\ee
follows from direct calculations.

However, in general 3D cases, the topology of the cross section may rather complicated and one can hardly represent the formula of the solution as \eqref{2dpoi}  or \eqref{Hagen}. Whether the gradient bound \eqref{E2} still holds for 3D Poiseuille flow in a sufficiently regular pipe is an interesting question. In this paper, we give a positive answer to this question. With the help of this, we prove the Poiseuille flow given in \eqref{HP1} is the unique solution to the problem \eqref{NS}--\eqref{NBC}--\eqref{Fluxc} provided $\Phi$ is no larger than a constant that is independent with the friction ratio $\al$.

Now we outline the idea of proving the main result of the current paper. First and the most important is to prove the upper bound \eqref{E1777} holds for Poiseuille flow in the general smooth pipe $\mD$. To do this, we first rescale \eqref{HP1} to get
\be\label{sys}
\left\{
\begin{aligned}
&-\Dl_h\varphi_\al=Constant_\al,\q\text{in }\Sigma\,,\\
&\f{\p\varphi_\al}{\p\bar{\bl{n}}}+\al\varphi_\al=0,\q\text{on }\p\Sigma\,,\\
&\int_{\Sigma}\varphi_\al dx_h=1\,.
\end{aligned}
\right.
\ee
To prove $\|\na_h\varphi_\al\|_{L^2(\Sigma)}\lesssim \al$ for $\al$ being sufficiently small, we expand $\varphi_\al$ as:
\be\label{Series}
\varphi_\al(x_h):=\sum_{n=0}^\i\al^n\phi_n(x_h)\,,
\ee
where $\phi_{0}=|\Sigma|^{-1}$, and $\phi_n$ satisfies
\[
\left\{
\begin{aligned}
&-\Dl_h\phi_{n}=Const_n,\q\text{in }\Sigma\,;\\
&\f{\p \phi_{n}}{\p\bar{\bl{n}}}=-\phi_{n-1},\q\text{on }\p\Sigma\,;\\
&\int_{\Sigma}\phi_{n}dx_h=0\,,
\end{aligned}
\right.\q\text{for}\q n\in\mathbb{Z}_+\,.
\]
The key observation of this decomposition is the sequence \eqref{Series} is convergent when $\al$ is small, and the zero-order term $\phi_0$ is a constant which makes no sense after taking the gradient. Meanwhile, the uniform boundedness of $\|\na_h\varphi_\al\|_{L^2(\Sigma)}$ for sufficient large $\al$ is concluded by an energy estimate of system \eqref{sys}. After showing the estimate
\[
\|\na{P}_\Phi\|_{L^2(\Sigma)}\leq\f{C\al\Phi}{1+\al}\,,
\]
the uniqueness is proved by deriving a Saint-Venant type estimate that was first introduced by Ladyzhenskaya and Solonnikov in \cite{LaSo}.

%Finally, our approach in the current paper fails for the general 3D Leray's problem subject to the Navier-slip boundary condition. The difficulty occurs on the compact distorted part in the middle of the pipe. We have some detailed discussions in Section \ref{Diss}, and it seems that this expose the essential difference between the 2D and higher dimensional flow subject to the Navier-slip boundary condition.

The rest of this paper is organized as follows. Section \ref{PREL} contains the preliminary work of the proof, in which some useful lemmas are presented. Section \ref{SEC3} is devoted to deriving the key estimate of the Poiseuille flow \eqref{HP1}. In Section \ref{SEC4}, we finish the proof of the uniqueness of the solution, and prove the critical flux $\Phi_0$ in no smaller than $\f{\pi}{16}$ if the cross section of the pipe is a unit disk. %Finally, some discussions on why the method given in the current paper fails for the general Leray's problem with a distorted part in the middle of the pipe are presented in Section \ref{Diss}.

\subsection{Notations}
We finish this section with a list of notations that will appear throughout the paper.
\begin{itemize}
\item $C_{a,b,...}$ denotes a positive constant depending on $a,\,b,\,...$ which may be different from line to line. $A\lesssim B$ means $A\leq CB$, and $A\lesssim_{a,b,...}B$ denotes $A\leq C_{a,b,...}B$. Meanwhile, $A\simeq B$ means both $A\lesssim B$ and $B\lesssim A$, in other words, $A$ and $B$ can be controlled by each other.

\item In the standard Euclidean coordinates framework $\boldsymbol{e_1},\ \boldsymbol{e_2}$ and $\boldsymbol{e_3}$, for a 3D vector field $\bl{w}$, we denote $\bl{w}=w_1\boldsymbol{e_1}+w_2\boldsymbol{e_2}+w_3\boldsymbol{e_3}$.

\item We set $\boldsymbol{n}=(n_1,n_2,0)$ is the unit outer normal vector on $\p\mD=\p\Sigma\times\mR$, where $\bar{\boldsymbol{n}}=(n_1,n_2)$ is the unit outer normal vector on $\p\Sigma$. Meanwhile, $\bl{\tau}=(\tau_1,\tau_2,0)$ is the unit tangent vector on $\p\Sigma\times\mR$ that is orthogonal to $\bl{e_3}$, and $(\bl{\tau},\bl{e_3},\bl{n})$ form a right-hand coordinate system. See Figure \ref{FIG1}.

\item $\mathfrak{H}$ stands for a multi-index such that $\mathfrak{H}=(h_1,h_2,h_3)$ where $h_1,h_2,h_3\in\mathbb{N}\cup\{0\}$ and $|\mathfrak{H}|=h_1+h_2+h_3$, $\nabla^\mathfrak{H}=\p_{x_1}^{h_1}\p_{x_2}^{h_2}\p_{x_3}^{h_3}$.

\item We use standard notations for Lebesgue and Sobolev functional spaces in $\mathbb{R}^3$: For $1\leq p\leq\infty$ and $k\in\mathbb{N}$, $L^p$ denotes the Lebesgue space with norm
\[
\|f\|_{L^p}:=
\lt\{
\begin{aligned}
&\left(\int_{\mathbb{R}^3}|f(x)|^pdx\right)^{1/p},\quad 1\leq p<\infty,\\
&\mathop{ess sup}_{x\in\mathbb{R}^3}|f(x)|,\quad\quad\quad\quad p=\infty.\\
\end{aligned}
\rt.
\]

\item $H^m$ denotes the $L^2$-based Sobolev space with its norm
\[
\begin{split}
\|f\|_{H^m}:=&\sum_{0\leq|\mathfrak{H}|\leq m}\|\nabla^\mathfrak{H} f\|_{L^2}\,.\\
\end{split}
\]

\item For any $\zeta>1$, we define
$$
\mathcal{D}_\zeta:=\left\{x \in \mathcal{D}:-\zeta<x_3<\zeta\right\}
$$
the truncated pipe with the length equals to $2\zeta$. Meanwhile, $\O^{\pm}_\zeta$ are defined by
\[
\Omega_\zeta^{+}:=\left(\mathcal{D}_\zeta-\mathcal{D}_{\zeta-1}\right) \cap\left\{x \in \mathcal{D}: x_3>0\right\}, \quad \Omega_\zeta^{-}:=\left(\mathcal{D}_\zeta-\mathcal{D}_{\zeta-1}\right) \cap\left\{x \in \mathcal{D}: x_3<0\right\}\,,
\]
respectively.
\end{itemize}

\section{Preliminaries}\label{PREL}
The following two lemmas are Poincar\'e inequalities of a vector $\bl{v}$ with only the normal part vanishes on the boundary. Detailed proof could be found in \cite[Lemma 2.5]{LPY2024SCM}.
\begin{lemma}\label{P2}
{\sl Let $\bl{v}=v_1\boldsymbol{e_1}+v_2\boldsymbol{e_2}+v_3\boldsymbol{e_3}$ be a $H^1$ vector field in $\Sigma\times I$, where $I\subset\mR$ is an interval. If $\bl{v}$ is divergence-free which satisfies $\bl{u}\cdot\boldsymbol{n}=0$ and
\[
\int_{\Sigma}v_3(x_h,x_3)dx_h=0\,,
\]
then we have
\[
\|\bl{v}\|_{L^2(\Sigma\times I)}\leq C\|\nabla_h \bl{v}\|_{L^2(\Sigma\times I)},
\]
where $C$ is an absolute positive constant.}
\end{lemma}

\qed

To estimate the pressure field in the Navier-Stokes equations, we need the following lemma for the divergence equation $\nabla\cdot \bl{V}=f$ in a truncated pipe. The detailed proof could be found in \cite{Bme1, Bme2}, also \cite[Chapter III]{Galdi2011}.

\begin{lemma}\label{LEM2.1}
{\sl Let $D=\Sigma\times [0,1]$, $f\in L^2(D)$ with
\[
\int_D fdx=0,
\]
then there exists a vector valued function $\bl{V}:\,D\to\mathbb{R}^3$ belongs to $H^1_0(D)$ such that
\be\label{LEM2.11}
\nabla\cdot \bl{V}=f,\q\text{and}\q\|\nabla \bl{V}\|_{L^2(D)}\leq C\|f\|_{L^2(D)}.
\ee
Here $C>0$ is an absolute constant.}
\end{lemma}

\qed

The following asymptotic estimate of a function that satisfies an ordinary differential inequality will be applied at the end of the proof. It was originally derived by Ladyzhenskaya and Solonnikov \cite{LaSo}. We also refer readers to \cite[Lemma 2.7]{LPY2024SCM} for a proof written in a relatively recent format.
\begin{lemma}\label{LEM2.3}
{\sl Let $Y(\zeta) \not \equiv 0$ be a nondecreasing nonnegative differentiable function satisfying
$$
Y(\zeta) \leqslant \Psi\left(Y^{\prime}(\zeta)\right), \quad \forall \zeta>0 .
$$

Here, $\Psi:[0, \infty) \rightarrow[0, \infty)$ is a monotonically increasing function with $\Psi(0)=0$, and there exist $C, \tau_1>0$ and $m>1$ such that
$$
\Psi(\tau) \leqslant C \tau^m, \quad \forall \tau>\tau_1 .
$$

Then
$$
\liminf _{\zeta \rightarrow+\infty} \zeta^{-\frac{m}{m-1}} Y(\zeta)>0\,.
$$}
\end{lemma}

\qed

\section{Asymptotic property of the Poiseuille flow with respect to $\al$}\label{SEC3}

In this section, we consider the Poiseuille flow \eqref{HP1} in $\mD=\Sigma\times\mR$. We show that the gradient bound for the Hagen-Poiseuille flow \eqref{E1777} remains valid for this general case. Here is the proposition:

\begin{proposition}\label{PROP2.2}
{\sl Given $\Phi\geq0$ be the flux, and let $\al\geq 0$ be the friction ratio in \eqref{NBC}. The Poiseuille flow defined in \eqref{HP1} satisfies the following estimate:
\be\label{EM}
\|\na{P}_\Phi\|_{L^2(\Sigma)}\leq\f{C\al\Phi}{1+\al}\,.
\ee
Here $C>0$ is a uniform constant which is independent with $\Phi$ or $\al$.}
\end{proposition}

\pf The existence and uniqueness of the scaler function $P_\Phi$ in \eqref{HP1} could be derived by routine methods of elliptic equations in a bounded smooth domain with the Robin boundary condition, here we omit the details. We only prove \eqref{EM} by showing the existence of $\tilde{C}_\Sigma>0$, such that
\be\label{PES}
\left\{
\begin{aligned}
\|\na{P}_\Phi\|_{L^2(\Sigma)}\leq{C\al\Phi},\q&\text{ for}\q 0<\al<\tilde{C}_\Sigma\,;\\[1mm]
\|\na{P}_\Phi\|_{L^2(\Sigma)}\leq{C\Phi},\q&\text{ for}\q \al\geq\tilde{C}_\Sigma\,.
\end{aligned}
\right.
\ee

We assume $\Phi>0$ since the remaining case is trivial. For simplicity, we denote $\varphi_\al:=\f{P_{\Phi}}{\Phi}$. This implies
\be\label{Evphi}
\left\{
\begin{aligned}
&-\Dl_h\varphi_\al=Constant_\al,\q\text{in }\Sigma\,;\\
&\f{\p\varphi_\al}{\p\bar{\bl{n}}}+\al\varphi_\al=0,\q\text{on }\p\Sigma\,;\\
&\int_{\Sigma}\varphi_\al dx_h=1\,.
\end{aligned}
\right.
\ee

To prove \eqref{PES}$_1$, we formally split
\be\label{ExS}
\varphi_\al(x_h):=\sum_{n=0}^\i\al^n\phi_n(x_h)\,,
\ee
where $\phi_{0}=|\Sigma|^{-1}$ that satisfies
\be\label{Ephi0}
\left\{
\begin{aligned}
&-\Dl_h\phi_{0}=Const_0,\q\text{in }\Sigma\,;\\
&\f{\p \phi_{0}}{\p\bar{\bl{n}}}=0,\q\text{on }\p\Sigma\,;\\
&\int_{\Sigma}\phi_{0}dx_h=1\,.
\end{aligned}
\right.
\ee
Meanwhile, for $n=1,2,3,...$, the function $\phi_n$ satisfies
\be\label{Ephin}
\left\{
\begin{aligned}
&-\Dl_h\phi_{n}=Const_n,\q\text{in }\Sigma\,;\\
&\f{\p \phi_{n}}{\p\bar{\bl{n}}}=-\phi_{n-1},\q\text{on }\p\Sigma\,;\\
&\int_{\Sigma}\phi_{n}dx_h=0\,.
\end{aligned}
\right.
\ee
Formally, one can verify $\varphi_\al$ satisfies \eqref{Evphi} by summing over \eqref{Ephi0} and \eqref{Ephin} for $n\in\mathbb{N}$. Since the solvability of \eqref{Ephin} is indeed classical, it remains to show the convergence of series \eqref{ExS} provided that $\al$ is sufficiently small. Integrating \eqref{Ephin}$_1$ over $\Sigma$ and using the Neumann-type boundary condition \eqref{Ephin}$_2$, one derives
\[
|Const_n|\cdot|\Sigma|=\Big|\int_{\Sigma}\Dl_h\phi_{n}dx_h\Big|=\Big|\int_{\p\Sigma}\f{\p\phi_{n}}{\p\bar{\bl{n}}}dS_h\Big|=\Big|\int_{\p\Sigma}\phi_{n-1}dS_h\Big|\leq|\p\Sigma|^{1/2}\|\phi_{n-1}\|_{L^2(\p\Sigma)}\,.
\]
Using the trace theorem (see e.g. \cite[Theorem 1 in Section 5.5]{Evans}), one arrives that
\be\label{EIcn}
|Const_n|\leq \tilde{C}_\Sigma\|\phi_{n-1}\|_{H^1(\Sigma)}\,,
\ee
which indicates
\be\label{Ecn}
|Const_n|\leq \tilde{C}^n_\Sigma
\ee
by an induction with \eqref{EIcn}.  Now multiplying $\phi_n$ on both sides of \eqref{Ephin}$_1$ and integration over $\Sigma$, one deduces
\[
-\int_{\Sigma}|\na_h\phi_n|^2dx_h=\int_{\p\Sigma}\phi_n\phi_{n-1}dS_h\,.
\]
Here we have applied the integration by parts and \eqref{Ephin}$_{2,3}$. This, together with the Poincar\'e inequality, the Young inequality and the trace theorem, indicates that
\[
\|\phi_n\|_{H^1(\Sigma)}^2\leq C_\Sigma\int_{\Sigma}|\na_h\phi_n|^2dx_h\leq C_\Sigma\|\phi_n\|_{L^2(\p\Sigma)}\|\phi_{n-1}\|_{L^2(\p\Sigma)}\leq\f{1}{2}\|\phi_n\|_{H^1(\Sigma)}^2+\f{\tilde{C}_\Sigma^2}{2}\|\phi_{n-1}\|_{H^1(\Sigma)}^2\,,
\]
which indicates
\be\label{Ephinh1}
\|\phi_n\|_{H^1(\Sigma)}\leq \tilde{C}_\Sigma\|\phi_{n-1}\|_{H^1(\Sigma)}\leq\,...\,\leq \tilde{C}_\Sigma^n\,.
\ee
Without loss of generality, we may assume constants $\tilde{C}_\Sigma$ on the far right of \eqref{Ecn} and \eqref{Ephinh1} are identical. Therefore, recalling \eqref{ExS}--\eqref{Ephin} and \eqref{Ecn}, if $\al<\tilde{C}_{\Sigma}^{-1}$, one has $\varphi_\al$ satisfies
\[
\left\{
\begin{aligned}
&-\Dl_h\varphi_\al=\sum_{n=0}^\i\al^nConst_n:=Constant_\al<\i,\q\text{in }\Sigma\,;\\
&\f{\p\varphi_\al}{\p\bar{\bl{n}}}+\al\varphi_\al=0,\q\text{on }\p\Sigma\,;\\
&\int_{\Sigma}\varphi_\al dx_h=1\,.
\end{aligned}
\right.
\]
And by \eqref{Ephinh1}, one has
\[
\|\varphi_\al\|_{H^1(\Sigma)}\leq\sum_{n=0}^\i\al^n\|\phi_n\|_{H^1}\leq\sum_{n=0}^\i(\al \tilde{C}_\Sigma)^n<\i\,.
\]
This indicates the validity of series \eqref{ExS}. Noticing that $\phi_0=|\Sigma|^{-1}$ is a constant, one concludes
\[
\|\na_h\varphi_\al\|_{L^2(\Sigma)}\leq\sum_{n=1}^\i\al^n\|\na_h\phi_n\|_{L^2(\Sigma)}\lesssim_\Sigma \al\,,
\]
which proves the case \eqref{PES}$_1$.

Now it remains to show $\|\na\varphi_\al\|_{L^2(\Sigma)}\lesssim_\Sigma 1$ for any $\al>\tilde{C}_\Sigma$. In fact, one only needs to show
\be\label{bdde}
\|\na\varphi_\al\|_{L^2(\Sigma)}\lesssim_\Sigma 1
\ee
for sufficient large $\al$ since the intermediate case is concluded by classical elliptic theory (see e.g. \cite{GT1998}). Thus \eqref{bdde} can be verified by showing
\be\label{Conv}
\na_h\varphi_\al\to\na_h\varphi_\i,\q\text{in }L^2(\Sigma),\q\text{as }\al\to +\i\,,
\ee
where $\varphi_\i$ satisfies
\[
\left\{
\begin{aligned}
&-\Dl_h\varphi_\i=Constant_\i,\q\text{in }\Sigma\,;\\
&\varphi_\i=0,\q\text{on }\p\Sigma\,;\\
&\int_{\Sigma}\varphi_\i dx_h=1\,.
\end{aligned}
\right.
\]
Denoting $\eta_\al:=\varphi_\al-\varphi_\i$, one derives
\[
\begin{split}
0&=(Constant_\al-Constant_\i)\int_{\Sigma}\eta_\al dx_h=-\int_{\Sigma}\eta_\al\Dl_h\eta_\al dx_h\\
&=\int_{\Sigma}|\na_h\eta_\al|^2dx_h+\al\int_{\p\Sigma}\eta_\al^2dS_h+\int_{\p\Sigma}\eta_\al\f{\p\varphi_\i}{\p\bar{\bl{n}}} dS_h\,.
\end{split}
\]
Using the Young inequality, this indicates
\[
\int_{\Sigma}|\na_h\eta_\al|^2dx_h+\f{\al}{2}\int_{\p\Sigma}\eta_\al^2dS_h\leq\f{1}{2\al}\int_{\p\Sigma}\big(\f{\p\varphi_\i}{\p\bar{\bl{n}}}\big)^2dS_h\,,
\]
which concludes \eqref{Conv}. This completes the proof of the proposition.

\qed

\section{Uniqueness}\label{SEC4}
In this section, we give the proof of the main result: the uniqueness of the Poiseuille flow provided the flux $\Phi$ is smaller than a given constant which is independent with the friction ration $\al$.

\subsection{Estimate of the pressure difference}
Below, the proposition shows an $L^2$ estimate related to the pressure in the truncated pipe $\O_Z^{+}$ or $\O_Z^{-}$ could be bounded by $L^2$ norm of $\nabla \bl{u}$.
\begin{proposition}\label{P2.4}
{\sl Let $({\bl{v}},{\tilde{\pi}})$ be an alternative weak solution of \eqref{NS} in the pipe $\mD$, subject to the Navier-slip boundary condition \eqref{NBC}. If the total flux
\[
\int_{\mD\cap\{x_3=z\}}{\bl{v}}(x_h,z)\cdot\boldsymbol{e_3}dx_h=\Phi,\q\text{for any } z\in\mR,
\]
then the following estimate the pressure difference holds:
\[
\left|\int_{\O^{\pm}_K}\big({\tilde{\pi}}-\pi_P\big)w_3dx\right|\leq C_\mD\left(\|\bl{P}_\Phi\|_{L^4(\O^{\pm}_K)}\|\nabla \bl{w}\|^2_{L^2(\O^{\pm}_K)}+\|\nabla \bl{w}\|_{L^2(\O^{\pm}_K)}^2+\|\nabla \bl{w}\|_{L^2(\O^{\pm}_K)}^3\right),\forall K\in\mR,
\]
where $\bl{w}:={\bl{v}}-{\bl{P}_\Phi}$, and $C_\mD>0$ is a constant independent of $K$.}
\end{proposition}
\pf In the following, the upper index ``$\pm$" of the domain and the tilde above the pressure will be canceled for simplicity. Noticing
\[
\int_{\mD\cap\{x_3=z\}}w_3(x_h,z)dx_h\equiv0,\q\forall z\in\mR,
\]
integrating over $z\in[K-1,K]$, we deduce that
\[
\int_{\O_K}w_3dx=0,\q\forall K>1.
\]
Using Lemma \ref{LEM2.1}, one derives that there exists a vector field $\bl{V}$ satisfying \eqref{LEM2.11} with $f=w_3$. By the stationary Navier-Stokes equations, one derives
\[
\int_{\O_K}\big({\pi}-\pi_P\big)w_3dx=-\int_{\O_K}\nabla \big({\pi}-\pi_P\big)\cdot \bl{V}dx=\int_{\O_K}\left(\bl{w}\cdot\na \bl{w}+\bl{P}_\Phi\cdot\nabla \bl{w}+\bl{w}\cdot\nabla \bl{P}_\Phi-\Dl \bl{w}\right)\cdot \bl{V}dx.
\]
This indicates
\[
\int_{\O_K}\big({\pi}-\pi_P\big)w_3dx=\sum_{i,j=1}^3\int_{\O_K}(\p_{x_i}w_j-w_iw_j)\p_{x_i}V_jdx-\sum_{i=1}^3\int_{\O_K}P_\Phi(w_i\p_{x_3}V_i+w_i\p_{x_i}V_3)dx.
\]
By applying H\"older's inequality and \eqref{LEM2.11} in Lemma \ref{LEM2.1}, one deduces that
\be\label{EP1}
\left|\int_{\O_K}\big({\pi}-\pi_P\big)w_3dx\right|\leq C\left(\|\na \bl{w}\|_{L^2(\O_K)}+\|\bl{w}\|_{L^4(\O_K)}^2+\|\bl{P}_\Phi\|_{L^\i(\O_K)}\|\bl{w}\|_{L^2(\O_K)}\right)\|w_3\|_{L^2(\O_K)}.
\ee
Since $w_3$ has a zero mean value on each cross-section $\Sigma$, and $(\bl{w}-w_3\bl{e_3})$ satisfies
\[
(\bl{w}-w_3\bl{e_3})\cdot\bl{n}=0,\quad\text{for any } x\in\p\mD\cap\p\O_K,
\]
the vector $\bl{w}$ enjoys the Poincar\'e inequality
\be\label{PON3}
\|\bl{w}\|_{L^2(\O_K)}\leq C_\mD  \|\na_h \bl{w}\|_{L^2(\O_K)}.
\ee
Applying \eqref{PON3} and the Gagliardo-Nirenberg interpolation, one concludes the following estimate from \eqref{EP1}:
\[
\left|\int_{\O_K}\big({\pi}-\pi_P\big)w_3dx\right|\leq C_\mD\left(\|\bl{P}_\Phi\|_{L^\i(\O_K)}\|\nabla \bl{w}\|^2_{L^2(\O_K)}+\|\nabla \bl{w}\|_{L^2(\O_K)}^2+\|\nabla \bl{w}\|_{L^2(\O_K)}^3\right).
\]

\qed

\subsection{Main estimates}
\q\ Noticing $\bl{P}_\Phi$ is an exact solution of \eqref{NS}, by subtracting the equation of $\bl{P}_\Phi$ from the equation of ${\bl{v}}$, one finds
\be\label{SUBT}
\bl{w}\cdot\nabla \bl{w}+\bl{P}_\Phi\cdot\nabla \bl{w}+\bl{w}\cdot\nabla \bl{P}_\Phi+\na\big({\pi}-\pi_P\big)-\Dl \bl{w}=0.
\ee
Multiplying $\bl{w}$ on both sides of \eqref{SUBT}, and integrating on $\mathcal{D}_\zeta$, one derives
\be\label{Maint0}
\int_{\mathcal{D}_\zeta}\bl{w}\cdot\Dl \bl{w}dx=\int_{\mathcal{D}_\zeta}\bl{w}\cdot\big(\bl{w}\cdot\nabla \bl{w}+\bl{P}_\Phi\cdot\nabla \bl{w}+\bl{w}\cdot\nabla \bl{P}_\Phi+\nabla\big({\pi}-\pi_P\big)\big)dx.
\ee
Using the divergence-free property and the Navier-slip boundary condition of $\bl{v}$ and ${\bl{P}_\Phi}$, one deduces
\be\label{Maint1}
\begin{split}
\int_{\mathcal{D}_\zeta}\bl{w}\cdot\Dl \bl{w}dx&=\int_{\mD_\zeta}w_i\p_{x_j}(\p_{x_j}w_i+\p_{x_i}w_j)dx\\
&=-\sum_{i,j=1}^3\int_{\mD_\zeta}\p_{x_j}w_i(\p_{x_j}w_i+\p_{x_i}w_j)dx+\sum_{i,j=1}^3\int_{\p\mD_\zeta}w_in_j(\p_{x_j}w_i+\p_{x_i}w_j)dx\\
&=-2\int_{\mD_\zeta}|\mathbb{S} \bl{w}|^2dx-\al\int_{\p\mD_\zeta\cap\p\mD}\left(|w_{\tau}|^2+|w_{3}|^2\right)dS\\
&\hskip 1cm+\sum_{i=1}^3\int_{\Sigma\times\{x_3=\zeta\}}w_i(\p_{x_3}w_i+\p_{x_i}w_3)dx_h-\sum_{i=1}^3\int_{\Sigma\times\{x_3=-\zeta\}}w_i(\p_{x_3}w_i+\p_{x_i}w_3)dx_h.
\end{split}
\ee
%which indicates
%\[
%\int_{\p^C\mD_Z}|u_\tau|^2dS=-\f{2}{\al}\int_{\mD_Z}|\mathbb{S}u|^2dx+\f{1}{\al}\left(\int_{\Sigma_Z}u_i(\p_3u_i+\p_iu_3)dx_h-\int_{\Sigma_{-Z}}u_i(\p_3u_i+\p_iu_3)dx_h\right)-\f{1}{\al}\int_{\mD_Z}u\Dl udx.
%\]
Here $\bl{n}=(n_1,n_2,n_3)$ is the unit outer normal vector on $\p\mD$.

On the other hand, using integration by parts, we may derive alternatively:
\be\label{Maint20}
\begin{split}
\int_{\mathcal{D}_\zeta}\bl{w}\cdot\Dl \bl{w}dx=-\int_{\mD_\zeta}|\nabla \bl{w}|^2dx+\un{\f{1}{2}\int_{\p\mD_\zeta}\nabla |\bl{w}|^2\cdot\boldsymbol{n}dS}_{I_1},
\end{split}
\ee
where
\[
\begin{split}
I_1=&\un{\f{1}{2}\int_{\p\mD_\zeta\cap\p\mD}\nabla |\bl{w}|^2\cdot\boldsymbol{n}dS}_{I_{2}}+\f{1}{2}\left(\int_{\mD\cap\{x_3=\zeta\}}\p_{x_3}|\bl{w}|^2(x_h,\zeta)dx_h-\int_{\mD\cap\{x_3=-\zeta\}}\p_{x_3}|\bl{w}|^2(x_h,-\zeta)dx_h\right).
\end{split}
\]
Using the global natural coordinates \eqref{NaCo}, one writes
\[
\nabla |\bl{w}|^2=\p_{\bl{\tau}}|\bl{w}|^2\boldsymbol{\tau}+\p_{x_3}|\bl{w}|^2\boldsymbol{e_3}+\p_{\bl{n}}|\bl{w}|^2\boldsymbol{n},\q\text{on }\p\mD_{\zeta}\cap\p{\mD}.
\]
Thus by \eqref{NBCM1}, one has $I_{2}$ satisfies
\be\label{T12}
|I_{2}|\leq\left|\int_{\p\mD_{\zeta}\cap\p{\mD}}\big(w_{\tau}\left(\al-\kappa(x)\right)w_{\tau}+\al|w_3|^2\big)dS\right|\leq \left(\al+\|\kappa\|_{L^\i(\p\mD)}\right)\int_{\p\mD_{\zeta}\cap\p{\mD}}|\bl{w}_{\mathrm{tan}}|^2dS.\\
\ee
Substituting \eqref{T12} in \eqref{Maint20}, one concludes that
\be\label{Maint2}
\int_{\mathcal{D}_\zeta}\bl{w}\cdot\Dl \bl{w}dx\leq-\int_{\mD_\zeta}|\nabla \bl{w}|^2dx+\left(\al+\|\kappa\|_{L^\i(\p\mD)}\right)\int_{\p\mD_\zeta\cap\p\mD}|\bl{w}_{\mathrm{tan}}|^2dS+C\int_{\mD\cap\{x_3=\pm\zeta\}}|\bl{w}||\nabla \bl{w}|dx_h.
\ee
Now we combine estimates \eqref{Maint1} and \eqref{Maint2} by calculating
\[
\eqref{Maint1}\times \left(\al+\|\kappa\|_{L^\i(\p\mD)}\right)+\eqref{Maint2}\times\al,
\]
and derive that
\be\label{Maint4}
\begin{split}
\left(2\al+\|\kappa\|_{L^\i(\p\mD)}\right)\int_{\mD_\zeta}\bl{w}\cdot\Dl \bl{w}dx\leq&-2\left(\al+\|\kappa\|_{L^\i(\p\mD)}\right)\int_{\mD_\zeta}|\mathbb{S}\bl{w}|^2dx-\al\int_{\mD_\zeta}|\nabla \bl{w}|^2dx\\
&+C_{\al,\kappa}\int_{\mD\cap\{x_3=\pm\zeta\}}|\bl{w}||\na \bl{w}|dx_h\,.
\end{split}
\ee

It remains to estimate the right-hand side of \eqref{Maint0}. Applying integration by parts, one derives
\be\label{Maint3}
\begin{split}
\int_{\mathcal{D}_\zeta}\bl{w}\cdot\big(\bl{w}\cdot\nabla \bl{w}+\nabla \big({\pi}-\pi_P\big)\big)dx=&\int_{\mD\cap\{x_3=\zeta\}}w_3\left(\f{1}{2}|\bl{w}|^2+\big({\pi}-\pi_P\big)\right)dx\\
&-\int_{\mD\cap\{x_3=-\zeta\}}w_3\left(\f{1}{2}|\bl{w}|^2+\big({\pi}-\pi_P\big)\right)dx.
\end{split}
\ee
Using integration by parts, one derives
\be\label{Maint7}
\int_{\mathcal{D}_\zeta}\bl{P}_\Phi\cdot\nabla \bl{w}\cdot \bl{w}dx=\f{1}{2}\int_{\mD_\zeta}\mathrm{div}(\bl{P}_\Phi|\bl{w}|^2)dx=\f{1}{2}\int_{\mD\cap\{x_3=\zeta\}}P_\Phi|\bl{w}|^2dx_h-\f{1}{2}\int_{\mD\cap\{x_3=-\zeta\}}P_\Phi|\bl{w}|^2dx_h.
\ee
Finally, applying H\"older's inequality, one derives
\[
\begin{split}
\left|\int_{\mathcal{D}_\zeta}\bl{w}\cdot\nabla \bl{P}_\Phi\cdot \bl{w}dx\right|=&\left|\int_{-\zeta}^\zeta\int_{\Sigma}\bl{w}_h\cdot\na_h P_\Phi w_3dx_hdx_3\right|\leq\|\na P_{\Phi}\|_{L^2(\Sigma)}\int_{-\zeta}^\zeta\|\bl{w}(\cdot,x_3)\|_{L^4(\Sigma)}^2dx_3\,.
\end{split}
\]
Using the result in Proposition \ref{PROP2.2} and the Gagliardo-Nirenberg interpolation, one derives
\be\label{Maint5}
\begin{split}
\left|\int_{\mathcal{D}_\zeta}\bl{w}\cdot\nabla \bl{P}_\Phi\cdot \bl{w}dx\right|&\leq\f{C\Phi\al}{1+\al}\int_{-\zeta}^\zeta\|\bl{w}(\cdot,x_3)\|_{L^2(\Sigma)}\left(\|\na_h\bl{w}(\cdot,x_3)\|_{L^2(\Sigma)}+\|\bl{w}(\cdot,x_3)\|_{L^2(\Sigma)}\right)dx_3\\
&\leq\f{C\Phi\al}{1+\al}\int_{\mD_\zeta}|\nabla \bl{w}|^2dx\,.
\end{split}
\ee
Here in the second inequality, we have applied Poincar\'e inequality Lemma \ref{P2} and the Cauchy-Schwarz inequality for integration on $x_3$ direction. Substituting \eqref{Maint4}, \eqref{Maint3}, \eqref{Maint7} and \eqref{Maint5} in \eqref{Maint0}, one derives
\[
\begin{split}
&\left(1-\f{C\left(2\al+\|\kappa\|_{L^\i(\p\mD)}\right)\Phi}{1+\al}\right)\int_{\mD_\zeta}|\nabla \bl{w}|^2dx\\
\leq& C_{\al,\kappa}\left(\int_{\mD\cap\{x_3=\pm\zeta\}}\big(|\bl{w}|(|\na \bl{w}|+|\bl{w}|^2)+|P_\Phi||\bl{w}|^2\big)dx_h\right.\\
&\left.-\int_{\mD\cap\{x_3=\zeta\}}w_3\big({\pi}-\pi_P\big)dx_h+\int_{\mD\cap\{x_3=-\zeta\}}w_3\big({\pi}-\pi_P\big)dx_h\right).
\end{split}
\]
Here the constant $C$ on the left is independent with $\al$ and $\Phi$. Now if $\Phi<<1$ being small enough such that $C\Phi\leq\f{1}{2}(2+\|\kappa\|_{L^\i(\p\mD)})^{-1}$, we have
\[
1-\f{C\left(2\al+\|\kappa\|_{L^\i(\p\mD)}\right)\Phi}{1+\al}\geq1-C\Phi\left(2+\|\kappa\|_{L^\i(\p\mD)}\right)\geq\f{1}{2}\,,
\]
which indicates
\[
\begin{split}
\int_{\mD_\zeta}|\nabla \bl{w}|^2dx\leq {C}_{\al,\kappa}&\left(\int_{\mD\cap\{x_3=\pm\zeta\}}\big(|\bl{w}|(|\na \bl{w}|+|\bl{w}|^2)+|P_\Phi||\bl{w}|^2\big)dx_h\right.\\
&\left.-\int_{\mD\cap\{x_3=\zeta\}}w_3\big({\pi}-\pi_P\big)dx_h+\int_{\mD\cap\{x_3=-\zeta\}}w_3\big({\pi}-\pi_P\big)dx_h\right).
\end{split}
\]
Integrating with $\zeta$ over $[K-1,K]$, one derives
\be\label{ET+0}
\int_{K-1}^K\int_{\mathcal{D}_\zeta}|\nabla \bl{w}|^2dxd\zeta\leq {C}_{\al,\kappa}\Bigg(\int_{\O_K^+\cup\O_K^-}\big(|\bl{w}|(|\na \bl{w}|+|\bl{w}|^2)+|P_\Phi||\bl{w}|^2\big)dx+\Big|\int_{\O_K^+\cup\O_K^-}w_3\big({\pi}-\pi_P\big)dx\Big|\Bigg).
\ee
Now we only handle integrations on $\O_K^+$ since the cases of $\O_K^-$ are similar. Using the Cauchy-Schwarz inequality and the Poincar\'e inequality Lemma \ref{P2}, one has
\be\label{ET+1}
\int_{\O_K^+}|\bl{w}||\na \bl{w}|dx\leq\|\bl{w}\|_{L^2(\O_K^+)}\|\nabla \bl{w}\|_{L^2(\O_K^+)}\leq C\|\nabla \bl{w}\|^2_{L^2(\O_K^+)}.
\ee
Moreover, by H\"older's inequality and the Gagliardo-Nirenberg inequality, one writes
\[
\int_{\O_K^+}|\bl{w}|^3dx\leq C_\mD\left(\|\bl{w}\|_{L^2(\O_K^+)}^{3/2}\|\na \bl{w}\|_{L^2(\O_K^+)}^{3/2}+\|\bl{w}\|_{L^2(\O_K^+)}^{3}\right),
\]
which follows by the Poincar\'e inequality that
%\be\label{ET+1.5}
\[
\int_{\O_K^+}|\bl{w}|^3dx\leq C\|\na \bl{w}\|_{L^2(\O_K^+)}^{3}.
\]
Meanwhile, by the Poincar\'e inequality and Proposition \ref{PROP2.2}, one deduces
\[
\int_{\O_K^+}|P_\Phi||\bl{w}|^2dx\leq\|P_\Phi\|_{L^\i}\|\bl{w}\|_{L^2(\O_K^+)}^2\leq C_{\al}\Phi\|\na \bl{w}\|_{L^2(\O_K^+)}^{2}\,.
\]
Recalling the estimate of pressure difference in Proposition \ref{P2.4}, one has
\be\label{ET+2}
\left|\int_{\O^{+}_K}\big({\pi}-\pi_P\big)w_3dx\right|\leq C_\mD\left(\|\bl{P}_\Phi\|_{L^\i(\O^{+}_K)}\|\nabla \bl{w}\|^2_{L^2(\O^{+}_K)}+\|\nabla \bl{w}\|_{L^2(\O^{+}_K)}^2+\|\nabla \bl{w}\|_{L^2(\O^{+}_K)}^3\right)\,.
\ee
For estimates \eqref{ET+1}--\eqref{ET+2} above, related inequalities on domain $\O_K^{-}$ could also be derived by the same approach. Substituting them in \eqref{ET+0}, one concludes
\be\label{FEST}
\int_{K-1}^K\int_{\mathcal{D}_\zeta}|\nabla \bl{w}|^2dxd\zeta\leq C_{\al,\Phi,\mD}\left(\|\nabla \bl{w}\|_{L^2(\O_K^+\cup\O_K^-)}^2+\|\nabla \bl{w}\|_{L^2(\O_K^+\cup\O_K^-)}^3\right).
\ee
\subsection{End of proof}

\q\ Finally, one concludes the uniqueness by the ordinary differential inequality Lemma \ref{LEM2.3}.  Defining
\[
Y(K):=\int_{K-1}^K\int_{\mathcal{D}_\zeta}|\nabla \bl{w}|^2dxd\zeta,
\]
\eqref{FEST} indicates
\[
Y(K)\leq C_{\al,\Phi,\mD}\left(Y'(K)+\left(Y'(K)\right)^{3/2}\right),\q\forall K\geq 1.
\]
Thus by Lemma \ref{LEM2.3}, we derive
\[
\liminf_{\zeta\to\infty}K^{-3}Y(K)>0,
\]
that is, there exists $C_0>0$ such that
\[
\int_{K-1}^K\int_{\mathcal{D}_\zeta}|\nabla \bl{w}|^2dxd\zeta\geq C_0K^3.
\]
However, this leads to a paradox with the condition \eqref{GROWC}. Thus, $Y(K)\equiv0$ for all $K\geq 1$, which proves $\bl{v}\equiv \bl{P}_\Phi$.

\qed

\subsection{An exact lower bound of the critical flux for cylindrical pipe}\label{Sec4.4}
In the case that $\mD$ is a unit cylindrical pipe, the solution of \eqref{HP1} is the generalized Hagen-Poiseuille flow given in \eqref{Hagen}. Here let us give an exact lower bound of critical flux $\Phi_0$ that guarantees the uniqueness. In this case, direct calculation shows
\be\label{Rem1}
|\na\bl{P}_{\text{disk},\Phi}(x)|\leq\f{4\al\Phi}{(\al+4)\pi}\,.
\ee
Indeed, the Poincar\'e inequality applied in \eqref{Maint5} consists of two parts: the $\bl{w}_h:=(w_1,w_2)$ part follows from Lemma \ref{P2} and the last component $w_3$ follows from the Poincar\'e inequality for functions with vanishing mean value. In the case of $\Sigma=D:=\{x_h\in\mR^2:\,|x_h|<1\}$, we first \textbf{claim} that
\be\label{EREM00}
\|\bl{w}_h\|^2_{L^2(D\times I)}\leq 2\left(\|\na_h\bl{w}_h\|^2_{L^2(D\times I)}+\|\p_{x_3}w_3\|^2_{L^2(D\times I)}\right)\,.
\ee
Here goes the reason. Integrating the identity
\[
\sum_{i, j=1}^2\left[\partial_{x_i}\left(w_ix_jw_j\right)-\partial_{x_i} w_ix_j w_j-|\boldsymbol{w}_h|^2-w_ix_j \partial_{x_i} w_j\right]=0
\]
on $x_h\in D$,  one deduces
\be\label{EREM0}
\int_{D}|\boldsymbol{w}_h|^2 d x_h=\sum_{i, j=1}^2 \int_{D} \partial_{x_i}\left(w_ix_jw_j\right) d x_h-\sum_{i, j=1}^2 \int_{D} \partial_{x_i} w_ix_jw_j d x_h-\sum_{i, j=1}^2 \int_{D} \partial_{x_i} w_jx_jw_id x_h\,.
\ee
The first term on the right hand side vanishes since $\bl{w}\cdot{\bl{n}}=0$ on $\p D$. Noticing $\mathrm{div}\,\bl{w}=0$ and using the Cauchy-Schwarz inequality, one has
\be\label{EREM1}
\Big|\sum_{i, j=1}^2\partial_{x_i} w_ix_jw_j\Big|\leq|\p_{x_3}w_3|\cdot|x_1w_1+x_2w_2|\leq|\p_3w_3||\bl{w}_h|\,,
\ee
and
\be\label{EREM2}
\Big|\sum_{i, j=1}^2\partial_{x_i} w_jx_jw_i\Big|\leq\Big(\sum_{i, j=1}^2|\p_{x_i}w_j|^2\Big)^{1/2}\Big(\sum_{i, j=1}^2|w_j|^2|x_i|^2\Big)^{1/2}\leq |\na_h\bl{w}_h||\bl{w}_h|\,.
\ee
Substituting \eqref{EREM1} and \eqref{EREM2} in \eqref{EREM0}, and applying the Young inequality, one deduces
\[
\int_{D}|\bl{w}_h|^2dx_h\leq\f{1}{2}\int_{D}|\bl{w}_h|^2dx_h+\left(\int_{D}|\na_h\bl{w}_h|^2dx_h+\int_{D}|\p_3w_3|^2dx_h\right)\,.
\]
Thus one concludes \eqref{EREM00} by integrating with $x_3$ variable over $I$. This proves the \textbf{claim}.

For the component $w_3$, using the Poincar\'e inequality with vanishing mean value (\cite[Section 3]{Payne1960}) and also integrating over $I$ for the third component, one arrives at
\be\label{EREM000}
\|w_3\|^2_{L^2(D\times I)}\leq \f{4}{\pi^2}\|\na_h w_3\|^2_{L^2(D\times I)}\,.
\ee
Thus by adding \eqref{EREM00} and \eqref{EREM000}, one concludes
\[
\|\bl{w}\|^2_{L^2(D\times I)}\leq2\|\na\bl{w}\|^2_{L^2(D\times I)}\,.
\]
Recalling \eqref{Rem1}, the estimate \eqref{Maint5} in the proof of Theorem \ref{th1} could be refined as
\[
\left|\int_{\mathcal{D}_\zeta}\bl{w}\cdot\nabla \bl{P}_{\text{disk},\Phi}\cdot \bl{w}dx\right|\leq\f{8\al\Phi}{(\al+4)\pi}\int_{\mD_\zeta}|\nabla \bl{w}|^2dx
\]
in the case $\Sigma=D$. According to the proof of Theorem \ref{th1} before, noticing that $\kappa\equiv1$ in this case, one only needs
\[
1-\f{8\left(2\al+1\right)\Phi}{(\al+4)\pi}>0
\]
to ensure the uniqueness of the generalized Hagen-Poiseuille flow. This indicates that the critical flux for generalized Hagen-Poiseuille flow satisfies
\[
\Phi_{\text{disk}, 0}\geq\f{\pi}{16}\thickapprox0.2\,.
\]
\section*{Acknowledgments}
\addcontentsline{toc}{section}{Acknowledgments}
\q Z. Li is supported by National Natural Science Foundation of China (No. 12001285) and Natural Science Foundation of Jiangsu Province (No. BK20200803).

 {\footnotesize

{\sc Z. Li: School of Mathematics and Statistics, Nanjing University of Information Science and Technology, Nanjing 210044, China, and Academy of Mathematics \& Systems Science, Chinese Academy of Sciences, Beijing 100190, China}

  {\it E-mail address:}  zijinli@nuist.edu.cn

\medskip
{
{\sc N. Liu: Academy of Mathematics \& Systems Science, Chinese Academy of Sciences, Beijing 100190, China}

  {\it E-mail address:}  liuning16@mails.ucas.ac.cn
}

\medskip

 {\sc T. Zhou: School of Mathematics and Statistics, Nanjing University of Information Science and Technology, Nanjing 210044, China}

  {\it E-mail address:}  zhoutaoran@nuist.edu.cn

}
\end{document}